\documentclass[12pt,draft]{amsart}

\usepackage{amssymb}

\begin{document}

\newtheoremstyle{hplain}%
  {\topsep}
  {\topsep}
  {\itshape}
  {}
  {\bfseries}
  {.}
  { }
  {\thmname{#2 }\thmnumber{#1}\thmnote{ \rm(#3)}}%
  
\newtheoremstyle{hdefinition}
  {\topsep}%
  {\topsep}%
  {\upshape}
  {}%
  {\bfseries}%
  {.}
  { }%
  {\thmname{#2 }\thmnumber{#1}\thmnote{ \rm(#3)}}%
 
\theoremstyle{hplain}
\newtheorem{thm}{Theorem}
\newtheorem{note}[thm]{Note}
\newtheorem{lemma}[thm]{Lemma}
\newtheorem{remark}[thm]{Remark}
\newtheorem{corollary}[thm]{Corollary}
\newtheorem{claim}{Claim}[thm]

\theoremstyle{hdefinition}
\newtheorem{definition}[thm]{Definition}
\newtheorem{question}[thm]{Question}

\newcommand{\restr}{\restriction}
\newcommand{\arr}{\longrightarrow}
\newcommand{\sub}{\subseteq}
\newcommand{\stick}{{}\sp {\bullet}\hskip -7.2pt \shortmid }
\newcommand{\squdia}{\square\hskip -13.5pt\diamondsuit}
\newcommand{\oom}{\stackrel{1-1}{\longrightarrow}}
\newcommand{\iso}{\stackrel{\sim}{\longrightarrow}}
\newcommand{\conc}{{}^\smallfrown}
\newcommand{\setm}{\setminus}
\newcommand{\Hx}[1]{\mathit{H}_{\omega_{#1}}}
\newcommand{\Hk}[1]{\mathit{H}_{#1}}
\newcommand{\Ht}{\mathit{H}_\theta}
\newcommand{\Hl}{\mathit{H}_\lambda}
\newcommand{\mbb}{\mathbb}
\newcommand{\mcal}{\mathcal}
\newcommand{\cc}{\twoheadrightarrow}
\newcommand{\notcc}{\twoheadrightarrow\!\!\!\!\!\!\!/\,\,\,\,\,\,\,}
\renewcommand{\frak}{\mathfrak}
\renewcommand{\theequation}{\thesection.\arabic{equation}}

\newcommand{\Htlk}{[\Ht]^{<\kappa}}

\numberwithin{equation}{section}

\newcommand{\rem}[1]{}
\newcommand{\comment}[1]{\marginpar{\scriptsize\raggedright #1}}
\renewcommand{\thesubsection}{\alph{subsection}}

\title[Chang's conjecture at supercompact cardinals]{Chang's
  conjecture may fail at supercompact cardinals}
\author{Bernhard K\"onig}
\address{\newline
  Universit\'e Paris 7\newline
  2 place Jussieu\newline
  75251 Paris Cedex 05\newline
  France}

\thanks{The author acknowledges a grant awarded by the French
  ministry of research.}


\begin{abstract}
  We prove a revised version of Laver's indestructibility theorem
  which slightly improves over the classical result. An application
  yields the consistency of
  $(\kappa^+,\kappa)\notcc(\aleph_1,\aleph_0)$ when $\kappa$ is
  supercompact. The actual proofs show that $\omega_1$-regressive
  Kurepa-trees are consistent above a supercompact cardinal even
  though ${\rm MM}$ destroys them on all regular cardinals. This
  rather paradoxical fact contradicts the common intuition.
\end{abstract}

\maketitle

\section{Introduction}

A structure $\frak{A}$ with a distinguished predicate $R$ is said to
be of {\em type $(\lambda,\kappa)$} if $|\frak{A}|=\lambda$ and
$|R^\frak{A}|=\kappa$. The relation
$(\lambda,\kappa)\rightarrow(\mu,\nu)$ then means that for every
structure of type $(\lambda,\kappa)$ there is an elementary equivalent
structure of type $(\mu,\nu)$. We have the classical theorem (see
\cite{chang90:_model}):

\begin{thm}[Vaught]\label{vaught-thm}
  $(\lambda^+,\lambda)\rightarrow(\aleph_1,\aleph_0)$ holds for all
  cardinals $\lambda$.
\end{thm}

Consider a variation of this notion: the principle
$(\lambda,\kappa)\cc(\mu,\nu)$ means that every structure of type
$(\lambda,\kappa)$ has an elementary substructure of type $(\mu,\nu)$.
The relation $(\aleph_2,\aleph_1)\cc(\aleph_1,\aleph_0)$ is usually
called {\em Chang's conjecture}.

\begin{thm}[Silver]\label{silver-thm}
  Chang's conjecture is independent of ${\rm ZFC}$.
\end{thm}

Silver's result (see \cite{kanamori78}) demonstrated in particular
that Chang's conjecture is related to large cardinals and
indiscernibles.  Its exact consistency strength was later established
by Donder \cite{donder79:_some} to be an $\omega_1$-Erd\H{o}s
cardinal. This left open the possibility of whether a version of
Chang's conjecture holds at or above a supercompact cardinal. The main
contribution of this paper is the following.

\begin{thm}\label{sup-cc-indep}
  $(\kappa^+,\kappa)\cc(\aleph_1,\aleph_0)$ is independent of ${\rm
    ZFC}$ even if $\kappa$ is a supercompact cardinal.
\end{thm}

The proof uses a technique developed by Laver known as the {\em
  indestructibility theorem} for supercompact
cardinals.\footnote{Recently, Apter and Hamkins have done a lot of
  work in this area and refined Laver's indestructibility theorem in
  various ways (see for example \cite{hamkins00}).} The key result
from \cite{laver78:_makin} reads as follows.

\begin{thm}\label{Lav-ind}
  A supercompact cardinal $\kappa$ can be made indestructible in the
  following sense: $\kappa$ is supercompact in any generic extension
  $V^\mbb{P}$ where $\mbb{P}$ is a $\kappa$-directed-closed partial
  ordering.
\end{thm}

But we will argue later in this article that indestructibility under
$\kappa$-directed-closed forcings is not enough for our purposes so we
have to establish a result slightly more general than Theorem
\ref{Lav-ind}. It was shown in \cite{koenig05:_kurep_namba} that a
supercompact cardinal is always destructible by $\lambda$-closed
forcing for arbitrarily large $\lambda$, so we cannot hope to
strengthen Laver's result by replacing the phrase
'$\kappa$-directed-closed' with '$\lambda$-closed' even if $\lambda$
is very large. The purpose of this note is to offer a strengthening of
Theorem \ref{Lav-ind} that goes into a different direction.

\begin{thm}\label{revised-intro}
  A supercompact cardinal $\kappa$ can be made indestructible in the
  following sense: $\kappa$ is supercompact in any generic extension
  $V^\mbb{P}$ where $\mbb{P}$ is a partial ordering that is 'almost
  everywhere' $\kappa$-directed-closed.
\end{thm}

The exact meaning of 'almost everywhere' will be made precise in
Definition \ref{def-F-directed}. Once established, Theorem
\ref{revised-intro} will be applied to prove Theorem
\ref{sup-cc-indep}. In a generalization of this argument, we later go
on to show that $(\lambda^+,\lambda)\cc(\aleph_1,\aleph_0)$ can fail
for all regular uncountable $\lambda$ even in the presence of a
supercompact cardinal.

\section{Notation}

The reader is assumed to have a strong background in set theory,
especially regarding Easton extensions. As general references we
recommend \cite{higherinf} and \cite{Kunen-intro}. Let us make some
remarks on the notations used in this paper.

We use an abbreviation in the context of elementary embeddings:
$j:M\arr N$ means that $j$ is a non-trivial elementary embedding from
$M$ into $N$ such that $M$ and $N$ are transitive. The {\em critical
  point} of such an embedding, i.e. the first ordinal moved by $j$, is
denoted by ${\rm cp}(j)$. We write $jx$ for $j(x)$ in a context where
too many parentheses might be confusing. An embedding $j:V\arr M$ is
called {\em $\lambda$-supercompact} if $\kappa={\rm cp}(j)$ is mapped
above $\lambda$ and $M$ is closed under $\lambda$-sequences or
equivalently if $j"\lambda\in M$. Remember that this is the same as
saying that there is an ultrafilter on $[\lambda]^{<\kappa}$ that is
{\em supercompact}, i.e. the set $\{x\in[\lambda]^{<\kappa}:\alpha\in
x\}$ is in the filter for every $\alpha<\lambda$ (fineness) and every
regressive function on a filter set is constant on a filter set
(normality). If $\kappa$ is an infinite cardinal, then the poset
$\mbb{P}$ is called {\em $\kappa$-directed-closed} if any directed
subset of size less than $\kappa$ has a lower bound in $\mbb{P}$ and
called {\em $\kappa$-closed} if any $\mbb{P}$-descending chain of
length less than $\kappa$ has a lower bound. $\mbb{P}$ is {\em
  strategically $\kappa$-closed} if Player Nonempty has a winning
strategy in the Banach-Mazur games of length less than $\kappa$. The
following important theorem from \cite{laver78:_makin} has become part
of the set-theoretic folklore.

\begin{thm}\label{lav-function}
  Let $\kappa$ be supercompact. Then there is an
  $f:\kappa\arr\Hk{\kappa}$ such that for every $x$ and every
  $\lambda\geq|{\rm TC}(x)|$ there is a $\lambda$-supercompact
  embedding $j:V\arr M$ such that $(jf)(\kappa)=x$.
\end{thm}

A function $f:\kappa\arr\Hk{\kappa}$ as in Theorem \ref{lav-function}
is usually called {\em Laver function} or {\em Laver diamond}. It is a
major tool in the proof of the Laver indestructibility theorem.

\section{Revised Indestructibility}

\begin{definition}\label{def-F-directed}
  If $\mbb{P}$ is a partial ordering then we always let
  $\theta=\theta_\mbb{P}$ be the least regular cardinal such that
  $\mbb{P}\in\Ht$. Say that an $X\in\Htlk$ is {\em $\mbb{P}$-complete}
  if every $(X,\mbb{P})$-generic filter has a lower bound in
  $\mbb{P}$. Define
  \begin{displaymath}
    \mcal{H}(\mbb{P})=\{X\in\Htlk: X\mbox{ is
    }\mbb{P}\mbox{-complete}\}.
  \end{displaymath}
  Then a partial ordering $\mbb{P}$ is called {\em almost
    $\kappa$-directed-closed} if $\mbb{P}$ is strategically
  $\kappa$-closed and $\mcal{H}(\mbb{P})$ is in every supercompact
  ultrafilter on $\Htlk$.
\end{definition}

Clearly, if a poset $\mbb{P}$ is $\kappa$-directed-closed then it is
almost $\kappa$-directed-closed. Thus, the following Theorem
\ref{rev-Lav-ind} is actually a bit stronger than the classical Laver
indestructibility. We will present applications later in which this
slight edge is crucial. Notice also that a closed unbounded subset of
$\Htlk$ is in every supercompact ultrafilter on $\Htlk$. But there are
more interesting examples:

\begin{lemma}\label{CF-claim}
  Let $\nu<\kappa\leq\theta$ be regular cardinals. Then $${\rm
    CF}(\geq\nu)=\{X\in\Htlk:{\rm cf}(\sup X\cap\kappa)\geq\nu\}$$ is
  in all supercompact ultrafilters on $\Htlk$. \hfill\qed
\end{lemma}

Lemma \ref{CF-claim} will be exploited later in the applications. Now
we prove the {\em revised indestructibility theorem}.

\begin{thm}[Revised Laver indestructibility]\label{rev-Lav-ind}
  A supercompact cardinal $\kappa$ can be made indestructible in the
  following sense: $\kappa$ is supercompact in any generic extension
  $V^\mbb{P}$ where $\mbb{P}$ is almost $\kappa$-directed-closed.
\end{thm}
\begin{proof}
  Of course, the proof is similar to the one given by Laver in
  \cite{laver78:_makin}. To show that our slight variation works, we
  give the whole proof in detail. Let $f:\kappa\arr\Hk{\kappa}$ be as
  in Theorem \ref{lav-function}.  We construct an Easton support
  iteration $\mbb{Q}_\kappa=(\mbb{Q}_\alpha:\alpha<\kappa)$ of length
  $\kappa$. During this iteration, we inductively define a poset
  $\mbb{Q}_\alpha$ and an ordinal $\lambda_\alpha$. If $\gamma$ is
  limit then $\mbb{Q}_\gamma$ will be the Easton support limit of
  $(\mbb{Q}_\alpha:\alpha<\gamma)$ and we define
  $\lambda_\gamma=\sup_{\alpha<\gamma}\lambda_\alpha$. In the
  successor step $\alpha+1$ we pick a $\mbb{Q}_\alpha$-name $\mbb{P}$
  for a partial ordering, where $\mbb{P}$ is trivial except when
  \begin{enumerate}
  \item[(i)] $\lambda_\xi<\alpha$ for all $\xi<\alpha$,
  \item[(ii)] $f(\alpha)=\langle\mbb{P},\lambda\rangle$, where
    $\mbb{P}$ is a $\mbb{Q}_\alpha$-name for a poset, and
  \item[(iii)] $\Vdash_{\mbb{Q}_\alpha}"\mbb{P}\mbox{ is almost
      $\alpha$-directed-closed.}"$
  \end{enumerate}
  If (i)-(iii) are satisfied then we let
  $\mbb{Q}_{\alpha+1}=\mbb{Q}_\alpha*\mbb{P}$ and
  $\lambda_{\alpha+1}=\lambda$. We claim that this works. The rest of
  the proof is designed to show that $\kappa$ is revised Laver
  indestructible in $V^{\mbb{Q}_\kappa}$.

  Now assume that $\mbb{P}$ is a $\mbb{Q}_\kappa$-name for a partial
  ordering that is almost $\kappa$-directed-closed, where $\theta$ is
  the least regular cardinal such that $\mbb{P}\in\Ht$. Remember that
  this means that $\mbb{P}$ is strategically $\kappa$-closed and that
  $\mcal{H}(\mbb{P})$ is in every supercompact ultrafilter on $\Htlk$.
  Let $\gamma\geq\kappa$. We need to find a supercompact ultrafilter
  on $[\gamma]^{<\kappa}$ in $V^{\mbb{Q}_\kappa*\mbb{P}}$. To this
  end, let $\lambda$ be a cardinal such that
  \begin{itemize}
  \item $\lambda>2^\theta$ and
  \item $\Vdash_{\mbb{Q}_\kappa*\mbb{P}}\lambda>
    2^{(\gamma^{<\kappa})}$.
  \end{itemize}
  Using the properties of the Laver function, pick in $V$ a
  $\lambda$-supercompact embedding $j:V\arr M$ such that
  $(jf)(\kappa)=\langle\mbb{P},\lambda\rangle$. Notice that
  $$j(\mbb{Q}_\alpha:\alpha\leq\kappa)=(\mbb{Q}_\alpha:\alpha\leq
  j\kappa).$$
  Now conditions (i)-(iii) are satisfied in $M$ when we
  replace $\kappa$ for $\alpha$ and $jf$ for $f$. Thus by
  elementarity, $$\mbb{Q}_{\kappa+1}=\mbb{Q}_\kappa*\mbb{P}$$
  and note
  that $\mbb{Q}_\delta$ is the trivial poset whenever
  $\kappa+1<\delta<\lambda$. So the final segment of the iteration
  $\mbb{Q}_{j\kappa}$ after the $(\kappa+1)$th step is by construction
  a strategically $\lambda$-closed forcing. But note that $j\mbb{P}$
  is also strategically $\lambda$-closed, so we can factor and get
  that $$\mbb{Q}_{j\kappa}* j\mbb{P}=\mbb{Q}_\kappa*\mbb{P}*\mbb{R}$$
  is such that the factor $\mbb{R}$ is strategically $\lambda$-closed.
  Now let $G\sub\mbb{P}$ be generic over the model
  $V^{\mbb{Q}_\kappa}$.
  \begin{claim}\label{G-extends-claim}
    $j"G$ extends to a condition $p_G$ in $j\mbb{P}$.
  \end{claim}
  \begin{proof}
    Notice that $G$ is an $(\Ht,\mbb{P})$-generic filter, so by
    elementarity we have that $j"G$ is a $(j"\Ht,j\mbb{P})$-generic
    filter. But
    \begin{displaymath}
      j"\Ht\in j\mcal{H}(\mbb{P})=\mcal{H}(j\mbb{P})
    \end{displaymath}
    by the assumption that $\mcal{H}(\mbb{P})$ is in every
    supercompact ultrafilter on $\Htlk$. This implies that $j"\Ht$ is
    $j\mbb{P}$-complete and therefore $j"G$ has a lower bound in
    $j\mbb{P}$.
  \end{proof}
  We have found a master condition $p_G$ for the final segment
  $\mbb{R}$ of the iteration. This gives rise to the next claim:
  \begin{claim}\label{ext-claim}
    The embedding $j:V\arr M$ can be extended to
    $$j:V^{\mbb{Q}_\kappa*\mbb{P}}\arr
    M^{\mbb{Q}_{j\kappa}*j\mbb{P}}.$$
  \end{claim}
  \begin{proof}
    This is the classical extension lemma of Silver. We only need to
    see that if $\tau$ is a $\mbb{Q}_\kappa*\mbb{P}$-name then $j\tau$
    becomes a $\mbb{Q}_{j\kappa}*j\mbb{P}$-name.
  \end{proof}
  The following is standard and we only sketch the argument: working
  in $V^{\mbb{Q}_\kappa*\mbb{P}}$ we construct a sequence
  $(r_\xi:\xi<\lambda)$ of $\mbb{R}$-conditions below $p_G$ such that
  for every $\mcal{X}\sub[\gamma]^{<\kappa}$ in
  $V^{\mbb{Q}_\kappa*\mbb{P}}$ there is an $r_\xi$ that decides if
  $j"\gamma$ is in $j\mcal{X}$ or not. Similarly, we decide the
  statements that guarantee normality of the following filter:
  $$\mcal{U}=\{\mcal{X}\sub[\gamma]^{<\kappa}:\exists\xi<\lambda \;\,
  r_\xi\Vdash j"\gamma\in j\mcal{X}\}.$$
  Then $\mcal{U}$ is a
  supercompact ultrafilter on $[\gamma]^{<\kappa}$ in
  $V^{\mbb{Q}_\kappa*\mbb{P}}$.
\end{proof}

\section{Regressive Kurepa-trees and
transfer principles}

Higher Kurepa-trees are natural counterexamples to model-theoretic
transfer principles. For example, it is easy to check that the
existence of an $\omega_2$-Kurepa-tree negates the relation
$(\aleph_3,\aleph_2)\cc(\aleph_2,\aleph_1)$. But in order to negate
the principle $(\aleph_3,\aleph_2)\cc(\aleph_1,\aleph_0)$, we need a
notion stronger than that of a regular $\omega_2$-Kurepa-tree. The
next definition is designed to serve this purpose.

\begin{definition}
  For any tree $T$ say that the level $T_\alpha$ is {\em
    non-stationary} if there is a function $f_\alpha:T_\alpha\arr
  T_{<\alpha}$ which is {\em regressive} in the sense that
  $f_\alpha(x)<_Tx$ for all $x\in T_\alpha$ and if $x,y\in T_\alpha$
  are distinct then $f_\alpha(x)$ or $f_\alpha(y)$ is strictly above
  the meet of $x$ and $y$.

  A $\lambda$-Kurepa-tree $T$ is called {\em $\gamma$-regressive} if
  $T_\alpha$ is non-stationary for all limit ordinals $\alpha<\lambda$
  with ${\rm cf}(\alpha)<\gamma$.
\end{definition}

The notion of a regressive Kurepa-tree was first introduced in
\cite{koenig05:_kurep_namba}. One can verify that an
$\omega_1$-regressive $\omega_2$-Kurepa-tree is a counterexample to
$(\aleph_3,\aleph_2)\cc(\aleph_1,\aleph_0)$. We actually prove the
more general:

\begin{lemma}\label{regKtree-notcc}
  Let $\kappa$ be regular and assume that there is an
  $\omega_1$-regressive $\kappa$-Kurepa-tree $T$. Then
  $(\kappa^+,\kappa)\notcc(\aleph_1,\aleph_0)$.
\end{lemma}
\begin{proof}
  Let $\mcal{B}$ be the set of cofinal branches of $T$ and consider
  the structure $(\mcal{B},T)$ which is of type $(\kappa^+,\kappa)$.
  Now assume towards a contradiction that
  $(\kappa^+,\kappa)\cc(\aleph_1,\aleph_0)$ would hold, so we find a
  substructure $$(\mcal{A},S)\prec(\mcal{B},T),$$
  where $\mcal{A}$ has
  size $\aleph_1$ and $S$ has size $\aleph_0$. Define
  $\delta=\sup({\rm ht}"S)$ and notice that ${\rm cf}(\delta)=\omega$.
  Hence $T_\delta$ is a non-stationary level of the tree $T$. A
  straightforward argument using the fact that there is a regressive
  1-1 function defined on $T_\delta$ shows that $\mcal{A}$ has size at
  most $|S|=\aleph_0$. This is a contradiction.
\end{proof}

Let us give a quick summary of what is known about regressive
Kurepa-trees.  A result from \cite{koenig05:_kurep_namba} says that
compact cardinals have considerable impact.

\begin{thm}\label{compact-noregKtree}
  Assume that $\kappa$ is a compact cardinal and $\lambda\geq\kappa$
  is regular. Then there are no $\kappa$-regressive
  $\lambda$-Kurepa-trees.
\end{thm}

Our goal is to show that a supercompact cardinal $\kappa$ is
consistent with the existence of an $\omega_1$-regressive
$\kappa$-Kurepa-tree. This would be in contrast to Theorem
\ref{compact-noregKtree}. But another theorem from
\cite{koenig05:_kurep_namba} indicates that we cannot succeed with the
classical Laver indestructibility.

\begin{thm}\label{MM-noregKtree}
  Under ${\rm MM}$, there are no $\omega_1$-regressive
  $\lambda$-Kurepa-trees for any uncountable regular $\lambda$.
\end{thm}

Indeed, if there were a $\kappa$-directed-closed partial ordering to
add an $\omega_1$-regressive $\kappa$-Kurepa-tree then such a forcing
would preserve ${\rm MM}$ in particular. But that contradicts Theorem
\ref{MM-noregKtree}. So we had to develop 'revised Laver
indestructibility' first in order to show the desired consistency.
This motivates the construction in the following section.

\section{An application}

\begin{lemma}\label{add-regKtree}
  Let $\nu\leq\lambda$ where $\lambda$ is regular. Then there is a
  $\lambda$-closed forcing $\mcal{K}^\lambda_\nu$ that adds a
  $\nu$-regressive $\lambda$-Kurepa-tree. Moreover, if $\kappa$ is
  supercompact and $\nu<\kappa\leq\lambda$ then $\mcal{K}^\lambda_\nu$
  is almost $\kappa$-directed-closed.
\end{lemma}
\begin{proof}
  We describe the forcing $\mcal{K}^\lambda_\nu$ and later show that it
  has the desired properties. One may assume the cardinal arithmetic
  $2^{<\lambda}=\lambda$, otherwise a preliminary Cohen-subset of
  $\lambda$ could be added. Conditions of $\mcal{K}^\lambda_\nu$ are
  pairs $(T,h)$, where
  \begin{enumerate}
  \item $T$ is a tree of height $\alpha+1$ for some $\alpha<\lambda$
    and each level has size $<\lambda$.
  \item $T$ is $\nu$-regressive, i.e. if $\xi\leq\alpha$ is of
    cofinality less than $\nu$ then $T_\xi$ is non-stationary over
    $T_{<\xi}$.
  \item $h:T_\alpha\arr\lambda^+$ is 1-1.
  \end{enumerate}

  The condition $(T,h)$ is stronger than $(S,g)$ if
  \begin{itemize}
  \item $S=T\restriction{\rm ht}(S)$.
  \item ${\rm rng}(g)\sub{\rm rng}(h)$.
  \item $g^{-1}(\nu)\leq_Th^{-1}(\nu)$ for all $\nu\in{\rm rng}(g)$.
  \end{itemize}

  A generic filter $G$ for $\mcal{K}^\lambda_\nu$ will produce a
  $\nu$-regressive $\lambda$-tree $T_G$ in the first coordinate and
  the sets $$b_\nu=\{x\in T_G:\mbox{there is }(T,h)\in G\mbox{ such
    that }h(x)=\nu\}$$
  for $\nu<\lambda^+$ form a collection of
  $\lambda^+$-many mutually different $\lambda$-branches through the
  tree $T_G$. Notice also that the standard arguments for
  $\lambda^+$-cc go through here as we assumed $2^{<\lambda}=\lambda$.

  It is shown in \cite{koenig05:_kurep_namba} that this forcing is
  $\lambda$-closed. So we are left with showing that if $\kappa$ is
  supercompact and $\nu<\kappa\leq\lambda$ then $\mcal{K}^\lambda_\nu$
  is almost $\kappa$-directed-closed. By Lemma \ref{CF-claim}, it
  suffices to show that whenever the set $X\in\Htlk$ is such that
  ${\rm cf}(\sup X\cap\kappa)\geq\nu$ then $X$ is $\mbb{P}$-complete.
  So assume that $X$ is like this and $K\sub\mbb{P}\cap X$ be an
  $(X,\mbb{P})$-generic filter. Define
  $$T_K=\bigcup\{T:\mbox{there is }(T,h)\mbox{ in }K\}$$ and
  $\delta={\rm ht}(T_K)=X\cap\kappa$. Then ${\rm cf}(\delta)\geq\nu$.
  Now extend $T_K$ by defining the level $T_\delta$ such that every
  $T_K$-branch colored by the filter $K$ has an extension in
  $T_\delta$. Note that the only problem here could be (2) of the
  definition of $\mcal{K}^\lambda_\nu$ but this condition does not
  apply to $T_\delta$ because ${\rm cf}(\delta)\geq\nu$.
\end{proof}

The supercompactness of $\kappa$ is not really used in the
construction of the forcing $\mcal{K}^\lambda_\nu$ but we included it
in the assumptions of Theorem \ref{add-regKtree} because the notion of
'almost $\kappa$-directed-closed', as we stated it, some\-how depends
on the supercompactness of $\kappa$. Using revised Laver
indestructibility, we get an immediate consequence.

\begin{thm}\label{Con-sc+regKtree}
  It is consistent with the supercompactness of $\kappa$ that there is
  an $\omega_1$-regressive $\kappa$-Kurepa-tree.
\end{thm}
\begin{proof}
  By Theorem \ref{rev-Lav-ind} and Lemma \ref{add-regKtree}.
\end{proof}

Theorem \ref{Con-sc+regKtree} is curious in the light of Theorem
\ref{MM-noregKtree}: even though ${\rm MM}$ forbids the existence of
$\omega_1$-regressive Kurepa-trees on all uncountable regular
cardinals, these trees are well consistent with a supercompact
cardinal. Imagine the following scenario: assume that $\kappa$ is
supercompact and $T$ an $\omega_1$-regressive $\kappa$-Kurepa-tree. If
we force ${\rm MM}$ over that model with the usual semiproper
iteration \cite{foreman88}, then $T$ becomes an $\omega_2$-Kurepa-tree
in the extension which would have non-stationary levels at all
ordinals that are $\omega$-cofinal in the ground model. But note that
the semiproper iteration to force ${\rm MM}$ will include Namba
forcing at many stages of this iteration, so the resulting model of
${\rm MM}$ contains far more $\omega$-cofinal ordinals than the ground
model in which $\kappa$ was still supercompact. This gives some
indication as to why the $\omega_1$-regressivity of $T$ is suddenly
impossible by Theorem \ref{MM-noregKtree}. Regarding this phenomenon,
see also \cite{koenig05:_kurep_namba}.

Going back to the transfer principles, our main concern was the
relation $(\kappa^+,\kappa)\cc(\aleph_1,\aleph_0)$ when $\kappa$ is
supercompact. A positive consistency result is straightforward and
well-known. Theorem \ref{measure-Coll} is probably folklore, but we
sketch the proof for convenience.

\begin{thm}\label{measure-Coll}
  Assume that $\kappa$ is regular and $\theta>\kappa$ is measurable.
  Then $$V^{{\rm Coll}(\kappa,<\theta)}\models
  (\kappa^+,\kappa)\cc(\aleph_1,\aleph_0).$$
\end{thm}
\begin{proof}
  Let $j:V\arr M$ be an embedding with critical point $\theta$. Using
  standard arguments, this embedding can be extended to
  \begin{displaymath}
    j:V^{{\rm Coll}(\kappa,<\theta)}\arr M^{{\rm
        Coll}(\kappa,<j\theta)}.
  \end{displaymath}
  It is now enough to show that, in $V^{{\rm Coll}(\kappa,<\theta)}$,
  every countable substructure $N$ can be {\em $\kappa$-end-extended},
  i.e. there is a proper elementary extension $\bar{N}$ such that
  $N\cap\kappa=\bar{N}\cap\kappa$. To see this, let $f_i\;(i<\omega)$
  enumerate all functions from $\theta$ to $\kappa$ that are in $N$.
  We may assume that $N$ contains everything in sight, so we have for
  all $i<\omega$,
  \begin{displaymath}
    (jf_i)(\theta)\in N\cap\kappa.
  \end{displaymath}
  By elementarity we can pick $\delta>\kappa$ such that
  $f_i(\delta)\in N\cap\kappa$ for all $i<\omega$. Now let $\bar{N}$
  be the Skolem Hull of $N\cup\{\delta\}$.
\end{proof}

\begin{corollary}\label{Con-sc+cc}
  The supercompactness of $\kappa$ is consistent with the
  model-theoretic transfer property
  $(\kappa^+,\kappa)\cc(\aleph_1,\aleph_0)$.
\end{corollary}
\begin{proof}
  Let $\kappa$ be a cardinal such that the supercompactness of
  $\kappa$ is indestructible by $\kappa$-directed-closed forcing and
  such that $\theta>\kappa$ is measurable. Then $\kappa$ remains
  supercompact in $V^{{\rm Coll}(\kappa,<\theta)}$. The corollary now
  follows from Theorem \ref{measure-Coll}.
\end{proof}

The following is the reverse to Corollary \ref{Con-sc+cc} and yields a
negative consistency result.

\begin{corollary}\label{Con-sc+ncc}
  The supercompactness of $\kappa$ is consistent with the
  model-theoretic transfer property
  $(\kappa^+,\kappa)\notcc(\aleph_1,\aleph_0)$.
\end{corollary}
\begin{proof}
  By Lemma \ref{regKtree-notcc} and Theorem \ref{Con-sc+regKtree}.
\end{proof}

\section{Global failure and GCH}

We generalize the technique of the last section to show that, in the
presence of a supercompact cardinal, the principle
$(\lambda^+,\lambda)\cc(\aleph_1,\aleph_0)$ can fail even for all
regular uncountable $\lambda$.

\begin{thm}\label{global-Ktree}
  It is consistent with the supercompactness of $\kappa$ that
  \begin{enumerate}
  \item ${\rm GCH}$ holds and
  \item there are $\omega_1$-regressive $\lambda$-Kurepa-trees for all
    regular $\lambda\geq\aleph_1$.
  \end{enumerate}
\end{thm}
\begin{proof}
  The reader is assumed to be familiar with the proof of Theorem
  \ref{rev-Lav-ind}. We only sketch the argument and content ourselves
  with pointing out the differences to the old construction. Start
  with a model in which $\kappa$ is supercompact and ${\rm GCH}$
  holds. Since we are only concerned with indestructibility for the
  forcings $\mcal{K}^\lambda_{\omega_1}$, our 'Laver preparation'
  consists only of posets of the same form. So define the Easton
  support iteration $\mbb{Q}=(\mbb{Q}_\alpha:\alpha\in{\rm Ord})$ by
  letting $\mbb{Q}_{\alpha+1}=\mbb{Q}_\alpha*\mbb{P}$, where
  \begin{itemize}
  \item $\mbb{P}=\mcal{K}^\alpha_{\omega_1}$ if $\alpha$ is regular
    uncountable and
  \item $\mbb{P}$ is trivial otherwise.
  \end{itemize}
  Using ${\rm GCH}$ and the standard Easton arguments, we know that
  all cardinals and cofinalities are preserved. The claim is that
  $\kappa$ is still supercompact in $V^\mbb{Q}$. Let
  $\gamma\geq\kappa$. Find a regular $\mu$ such that $\mbb{Q}$ factors
  into $\mbb{Q}_{\mu}*\mbb{R}$, where the final segment $\mbb{R}$ is
  $2^{(\gamma^{<\kappa})}$-closed and $\mu$ is much larger than
  $\gamma$. Now let $j:V\arr M$ be a $(2^\mu)^+$-supercompact
  embedding. Notice that $j(\mbb{Q}_{\mu})=\mbb{Q}_{j\mu}$ by
  construction. We can factor
  $$\mbb{Q}_{\mu}=\mbb{Q}_{\kappa}*\mbb{Q}_{[\kappa+1,\mu]}.$$
  If
  $H\sub\mbb{Q}_{\mu}$ is generic, we let $G$ be the restriction of
  $H$ to the final segment $\mbb{Q}_{[\kappa+1,\mu]}$. By the standard
  arguments, it is enough to show that $j"G$ extends to a condition
  $p_G$ in $j\mbb{Q}_{[\kappa+1,\mu]}=\mbb{Q}_{[j\kappa+1,j\mu]}$.
  This is similar to Claim \ref{G-extends-claim} and using the fact
  that $${\rm
    CF}(\geq\omega_1)\sub\mcal{H}(\mbb{Q}_{[\kappa+1,\mu]}).$$
  Note
  that $\mbb{Q}_{[\mu+1,j\mu]}/G$ is the same as
  $\mbb{Q}_{[\mu+1,j\mu]}$ below the condition $p_G$. Hence the
  closure properties of $\mbb{Q}_{[\mu+1,j\mu]}$ help us decide the
  properties of the supercompact ultrafilter on $[\gamma]^{<\kappa}$
  that lives in the model $V^{\mbb{Q}_{j\mu}}$ (compare with Claim
  \ref{ext-claim}). We have shown that $\kappa$ is
  $\gamma$-supercompact in $V^{\mbb{Q}_{\mu}}$. Since the final
  segment $\mbb{R}$ of the iteration is
  $2^{(\gamma^{<\kappa})}$-closed, it preserves
  $\gamma$-supercompactness and we are done.
\end{proof}

\begin{corollary}\label{global-notcc}
  It is consistent with the existence of a supercompact cardinal that
  $(\lambda^+,\lambda)\notcc(\aleph_1,\aleph_0)$ holds for all regular
  uncountable $\lambda$.
\end{corollary}
\begin{proof}
  By Lemma \ref{regKtree-notcc} and Theorem \ref{global-Ktree}.
\end{proof}

\bibliography{locco}\bibliographystyle{plain}

\end{document}